\documentclass[11pt]{article}
\usepackage[utf8]{inputenc}
\usepackage[english]{babel}
\usepackage{amsmath}
\usepackage{amsthm}

\usepackage{amsfonts}
\usepackage{amssymb}
\usepackage{makeidx}
\usepackage{graphicx}
\usepackage[left=2.6cm,right=2.6cm,top=3.5cm,bottom=2cm]{geometry}
\usepackage{enumerate}
\usepackage{mathrsfs}
\usepackage{multicol}

\newtheorem{theorem}{Theorem}[section]

\newtheorem{thmx}{Theorem}

\newtheorem{prop}[theorem]{Proposition}

\newtheorem{lemma}[theorem]{Lemma}

\newtheorem{remark}[theorem]{Remark}

\providecommand{\abs}[1]{\lvert#1\rvert}

\newcommand{\vsp}{\vspace{0,1 cm}}

\newcommand{\Z}{\mathbb Z}
\newcommand{\R}{\mathbb R}

\title{On the critical regularity of nilpotent groups acting\\ on the interval: the metabelian case} 

\begin{document}

\author{\textsc{Maximiliano Escayola and Crist\'obal Rivas }}

\maketitle

\begin{abstract}

Let $G$ be a torsion-free, finitely-generated, nilpotent and metabelian group. In this work we show that $G$ embeds into the group of orientation preserving $C^{1+\alpha}$-diffeomorphisms of the compact interval, for all $\alpha< 1/k$ where $k$ is the torsion-free rank of $G/A$ and $A$ is a maximal abelian subgroup. We show that in many situations the corresponding $1/k$ is  critical  in the sense that there is no embedding of $G$ with higher regularity. A particularly nice family where this happens, is the family of  $(2n+1)$-dimensional Heisenberg groups, for which we can show that the critical regularity  equals $1+1/n$.
\end{abstract}
\vspace{0.2cm}


\vspace{0.2cm}


\vspace{0.7cm}

\section{Introduction}

Given an integer $n\geq 0$ and a real number  $\alpha\in [0,1)$,  we denote by  $\text{Diff}_+^{\, n+\alpha}  ([0,1])$ the group of all orientation preserving $C^n$-diffeomorphisms of the closed interval $[0,1]$ whose $n$-th derivative is $\alpha$-H\"older continuous ($C^{n+\alpha}$-diffeomorphisms for short). Observe that with this notation the group $\text{Diff}_+^{\, 0}([0,1])$ is the group of all homeomorphisms of $[0,1]$ isotopic to the identity. Finally, observe that the family of groups $\text{Diff}_+^{\,n+\alpha}([0,1])$ is totally ordered by inclusion since $\text{Diff}_+^{\,n+\alpha}([0,1]) \supseteq \text{Diff}_+^{\,n'+\alpha'}([0,1]) $ if and only  if $n<n'$, or $n=n'$ and $\alpha\leq \alpha'$. 

We are interested in computing the critical regularity of an abstract group $G$ acting on the interval $[0,1]$. Recall that, given a group $G$, the {\em critical regularity} of  $G$ on $[0,1]$ is by definition 
$$\text{Crit}_{[0,1]}(G)=\text{sup}\{n+\alpha \mid n\geq 0,\alpha\in [0,1) \text{ and } G\text{ embeds into }\text{Diff}_+^{\,n+\alpha}([0,1])\},$$
where we set $\text{Crit}_{[0,1]}(G)=-\infty$ if $G$ does not embed into $\text{Diff}_+^{0}([0,1])$. The problem of computing the critical regularity of a group $G$ is quite natural and turns out to be very interesting in the case that $G$ is finitely-generated (the reader may wish to consult \cite{kk_book} for an introduction). For instance, we know from a theorem of Deroin, Kleptsyn, Navas \cite{DKN_acta} (see also \cite{dknp}) that every countable subgroup of $\text{Diff}^0_+([0,1])$ is conjugated to a group of bi-Lipschitz transformations, and hence $1\leq \text{Crit}_{[0,1]}(G)$ for every countable subgroup of $\text{Diff}^0_+([0,1])$ (for uncountable subgroups of $\text{Diff}^0_+([0,1])$ this is no longer true, see \cite{calegari}). However, the celebrated Stability Theorem of Thurston \cite{thurston} implies that every finitely-generated subgroup of $\text{Diff}^1_+([0,1])$ admits a surjective homomorphism onto the integers, and so not every group of homeomorphisms of the interval can be realized as a group of diffeomorphisms\footnote{Concrete examples of finitely-generated subgroups of $\text{Diff}_+^{\,0}([0,1])$ having trivial abelianization can be found in \cite{thurston,bergman, rivas-triestino}. However, Thuston's obstruction is not the only obstruction for $C^1$ smoothability as there are also known examples of finitely-generated and locally indicable groups having no {\em faithful} $C^1$ action on the interval, see \cite{calegari, navas-loc-ind,BMNR,kkr}.}.  Further obstructions appear in higher regularity: for $C^2$ there is the important Kopell obstruction \cite{kopell}, and between $C^1$ and $C^2$ there is the generalized Kopell obstruction from \cite{DKN_acta}.  In a related spirit, Kim and Koberda \cite{kk}, and later Mann and Wolff \cite{mann-wolff},  have shown that for every $n\geq 1$ and every $\alpha$ in $[0,1)$, there is a finitely-generated group whose critical regularity on $[0,1]$ is exactly $n+\alpha$.

In this work, we focus on actions on the interval of finitely-generated and torsion-free nilpotent groups (basic definitions will be recalled in \S \ref{sec on nilp metabelian}). Let $G$  be one such group. It follows from the work of Mal'cev that $G$ embeds into $\text{Diff}^0_+([0,1])$ (see, for instance, \cite[\S5.2]{robinson} and \cite[\S1.2]{GOD}), and we know from the work of Farb and Franks \cite{ff}, that every action of $G$ on $[0, 1]$ by homeomorphisms can be conjugated inside $\text{Diff}^1_+([0, 1])$ (see also the universal construction from Jorquera \cite{jorquera}). This was further refined by Parkhe \cite{p} who showed that actually any $C^0$-action of $G$ on $[0,1]$ can be conjugated inside $\text{Diff}^{1+\alpha}_+([0, 1])$ as long as $\alpha < 1/\tau$, where $\tau$ is the degree of the polynomial growth of the nilpotent group $G$. On the other hand, Plante and Thurston \cite{pt}, have shown that every nilpotent subgroup of $\text{Diff}^2_+([0, 1])$ must be abelian. So, if $G$ is a torsion-free, finitely-generated and nilpotent group which is non-abelian, then
$$1+1/\tau \leq  \text{Crit}_{[0,1]}(G) \leq 2.$$

The exact critical regularity of concrete nilpotent groups has been computed only in few cases and one important goal of this work is to provide new explicit computations of critical regularity for certain groups. Castro, Jorquera and Navas \cite{int4}, build a family of nilpotent abelian-by-cyclic groups whose critical regularity is 2. These examples can be made of arbitrarily large nilpotency degree, yet they are all metabelian (i.e. their commutator subgroup is abelian).  Jorquera, Navas and the second author showed in \cite{JNR}  that the critical regularity of $N_4$ (the group of $4$ by $4$ upper triangular matrices with $1$’s on the diagonal) is $1+1/2$.  We point out that at the time of this writing, $N_4$ is the only  torsion-free nilpotent group whose critical regularity is known and turns out not to be an integer. Note  that $N_4$ is also a metabelian group.

\vsp

The main purpose of this article is to exhibit many other nilpotent groups whose critical regularity is strictly between 1 and 2.  Our main technical result is an improvement of Parkhe's lower bound for the critical regularity in the class of finitely-generated, torsion-free nilpotent groups which are metabelian (see Remark \ref{rem parkhe}). For the statement, recall that the torsion-free rank of an abelian group $H$ is the dimension of the $\mathbb Q$-vector space $H\otimes \mathbb Q$. We denote this rank by $\text{rank}(H)$.

\begin{thmx} \label{teo lower bound} Let $G$ be a non-abelian, torsion-free, finitely-generated nilpotent group which is metabelian, and  let $A$ be a maximal abelian subgroup containing  $[G,G]$. If $k=\text{\normalfont rank}(G/A)$, then 
$$G \hspace{0,2cm}\text{embeds into}\hspace{0,2cm} \textup{Diff}_+^{1+\alpha}([0,1]),\hspace{0,2cm}\text{for all}\hspace{0,2cm}\alpha< 1/k.$$ 
In particular $1+ 1/k\leq \text{\normalfont Crit}_{[0,1]}(G)$.
\end{thmx}

\begin{remark}\label{rem parkhe} By the Bass-Guivarc'h formula \cite{bass, guivarch}, the degree of the polynomial growth of a nilpotent group $G$ is $\tau= \sum_{i\geq 1} i \text{ \normalfont rank}(\gamma_i/\gamma_{i+1}) $, where $G=\gamma_1\geqslant \gamma_2\geqslant ...$ is the lower central series of $G$. In particular, for a nilpotent group  $G$ as in Theorem \ref{teo lower bound} with maximal abelian subgroup $A$, we have that  $\text{\normalfont rank}(G/A)< \tau$. Hence the lower bound for $\text{\normalfont Crit}_{[0,1]}(G)$ in Theorem \ref{teo lower bound} is (strictly) greater than Parkhe's lower bound.

\end{remark}

The proof of Theorem \ref{teo lower bound} is given in Section \ref{sec lower bound}.  Taking inspiration from the abelian-by-ciclic action from \cite[\S 4]{int4}, in \S \ref{sec action by homeos} we build, for a metabelian and finitely-generated torsion-free nilpotent group $G$, a family of actions of $G$ on the interval $[0,1]$ by orientation preserving homeomorphisms. This is done by first building actions of $G$ on $\Z^{k+1}$ which preserve a lexicographic order and then  ``projecting'' them into the interval. In \S \ref{sec action by diffeos}, we use the Pixton-Tsuboi technique \cite{pixton, tsuboi} to show that these actions can be smoothed to actions by $C^{1+\alpha}$-diffeomorphisms for any $\alpha < 1/k$. This section closely follows the work \cite{int4}  and the main difference is that we don't have explicit polynomials in the construction of the actions, but only bounds on them (see Proposition {\ref{prop polynomial bound}). Although these actions may not be faithful, in \S \ref{sec faithful action} we explain how to glue some of these actions in order to obtain an  embedding of $G$ into $\textup{Diff}_+^{1+\alpha}([0,1])$ for any  $\alpha<1/k$.

\vsp

In some situations, even the lower bound in Theorem \ref{teo lower bound} is not sharp in the sense that there are groups for which the theorem applies yet their critical regularity is strictly greater than the predicted lower bound. This is related to the possibility of splitting  the group as a product of two groups each allowing an embedding with higher regularity. We provide an  easy example of this phenomenon in \S \ref{sec example with higher regularity}. However, in many cases we can ensure that the inequality in Theorem \ref{teo lower bound} is indeed optimal, and in \S4.1 and \S 4.2 we provide two families of examples where we can obtain upper bounds for the  regularity and hence compute the critical regularity.

The first family of examples are the $(2n+1)$-dimensional discrete Heisenberg groups, that we denote  $\mathscr{H}_{n}$. Recall that by definition $$\mathscr{H}_n:=\left\lbrace \begin{pmatrix}
1 & \vec{x} & c\\
\vec{0}^{\;t} & I_n & \vec{y}^{\; t}\\
0 & \vec{0} & 1
\end{pmatrix}\;: \vec{x},\vec{y}\in \mathbb{Z}^n\text{ and }c\in\mathbb{Z} \right\rbrace,$$
where $I_n$ is the identity matrix of size $n$ and $\vec{0}^{\;t}$, $\vec{y}^{\; t}$ are the transposes of $\vec{0}$, $\vec{y}$ respectively. It is easy to see that these groups are nilpotent of degree two and hence they are metabelian. Moreover, a maximal abelian subgroup $A$ of $\mathscr H_n$ is given by the set of matrices whose corresponding vector $\vec{x}=0$. In particular $\mathscr H_n/A$ has torsion-free rank equal to $n$.  For this family we  show in \S\ref{sec Heissemberg} that there is no embedding of $\mathscr H_n$ into $\text{Diff}_+^{1+\alpha}([0,1])$ for $\alpha>1/n$. In particular  we obtain

\begin{thmx}\label{teo Heissemberg} Let $\mathscr H_n$ be the $(2n+1)$-dimensional discrete Heisenberg group. Then
$$\textup{Crit}_{[0,1]}(\mathscr{H}_n)=1+\frac{1}{n}.$$
\end{thmx}

Finally,  in \S \ref{sec more examples} we produce examples of  metabelian and torsion-free nilpotent groups for which we can compute the critical regularity but whose nilpotency degree can be chosen to be arbitrarily large. More precisely we show

\begin{thmx}\label{teo clase nilp arbitraria}
For any integers $k$ and $d$ with $d>k$, there is a nilpotent group $G$ and a maximal abelian subgroup $A$ containing $[G,G]$ such that $d$ is the nilpotency degree of $G$,  $k$ is the torsion-free rank of $G/A$  and
$$\textup{Crit}_{[0,1]}(G)=1+\frac{1}{k}.$$

\end{thmx}

In both cases, the key to obtain an upper bound for the regularity is to use the internal algebraic structure of the groups in  order to be able to apply the generalized Kopell lemma from \cite{DKN_acta}.

\begin{remark} We know from  the results of Kim and Koberda \cite{kk} that for any real number $\alpha \geq 1$ there is a finitely-generated group whose critical regularity is exactly $\alpha$. However, in all known cases where the critical regularity of a torsion-free nilpotent group has been computed, it is of the form $1+1/n$ for some integer $n$. See Theorems B, C and  \cite{int4}, \cite{JNR}. So, we wonder if this is always the case for torsion-free and finitely-generated nilpotent groups (non necessarily metabelian). 
\end{remark}

\paragraph{Acknowledgements:} 
The first author acknowledges the hospitality of the Institut de  Mathématique de Bourgogne where part of this project was carried out. In particular he is grateful to Michele Triestino for fruitful discussion around this and other topics of mathematics. He is also grateful to Nicol\'as Matte Bon for discussions around the algebraic structure and actions of nilpotent groups.

The first author acknowledges the partial support of Project Gromeov ANR-19-CE40-0007 and Mathamsud 210020 ``Dynamical Groups Theory". Both authors are grateful to Andr\'es Navas for reading a preliminary version of this work. Both authors acknowledge the support of FONDECYT  1210155. Finally, both authors acknowledge the careful reading of the referees.

\section{Preliminaries on nilpotent groups and invariant orders}
\label{sec on nilp metabelian}

Given a group $G$ and two elements $f,g\in G$, we let  $[f,g]=fgf^{-1}g^{-1}$ denote the commutator of $f$ and $g$. Further, if $G$ is finitely-generated and $S$ is a finite generating set, an element of the form $[s_1,s_2]$ with $s_1,s_2\in S$ is called a simple commutator of weight 2. Inductively,  a \emph{simple commutator of weight} $n$ is defined  as an element of the form $$[s_1,\ldots,s_n]:=[s_1,[s_2,\ldots,s_n]],\hspace{0,4cm} s_1,\ldots,s_n\in S.$$ Note that given $n$, there exists only a finite number of simple commutators of weight $n$. 

Let $H$ and $K$ be subgroups of $G$. $[K,H]$ denotes the subgroup of $G$ generated by commutators $[g,h]$ with $g\in K$ and $h\in H$. The subgroup $[G,G]$ is called the {\em commutator subgroup} and we say that $G$ is {\em metabelian} if $[G,G]$ is abelian.  

Remember that the \emph{lower central series} of $G$ is $$G=\gamma_0\geqslant \gamma_1\geqslant \gamma_2\geqslant \cdots $$ where $\gamma_1=[G,G]$ and $\gamma_i=[G,\gamma_{i-1}]$; and the \emph{upper central series} of $G$ is
$$\{e\}= \zeta_0 \leqslant \zeta_1 \leqslant \zeta_2\leqslant \cdots  $$ 
where  $\zeta_{i}/\zeta_{i-1}=Z\left(G/\zeta_{i-1}\right)$, and $Z(G)$ denotes the center of $G$.

Let us recall some classic results about nilpotent groups. See \cite{robinson} for an in-depth exposition of them.   The group $G$ is \emph{nilpotent} of degree $n$ if $\zeta_n=G$ but $\zeta_{n-1}\neq G$. In this case it also happens that $\gamma_n=\{e\}$ but $\gamma_{n-1}\neq \{e\}$. Therefore, in a finitely-generated nilpotent group of degree $n$, we only have a finite number of simple commutators (for a fixed generating set). This is because all simple commutators of weight $n$ are trivial. 
 
 It is a result of Mal'cev that if $G$ is a torsion-free nilpotent group then the factors $\zeta_i/\zeta_{i-1}$ are also torsion-free for all $ i \in \{1,\ldots,n\}$ (see \cite[Proposition 5.2.19]{robinson}). 
Recall also that finitely-generated nilpotent groups are polycyclic, hence every subgroup of a finitely-generated nilpotent group is finitely-generated as well (see \cite[Proposition 5.4.12]{robinson}). Additionally nilpotent groups also satisfy that their non-trivial normal subgroups always intersect non-trivially the center of the group (see \cite[Proposition 5.2.1]{robinson}). An immediate consequence of this is the following useful result.

\begin{prop}\label{prop inyectiva}  Let $G$ be a nilpotent group, let $H$ be a group, and let $\varphi : G \rightarrow H$ be a group homomorphism. Then
$\varphi$ is injective if and only if $\varphi\mid_{Z(G)}$ (the restriction of $\varphi$ to $Z(G)$) is injective.
\end{prop}

For $g\in G$, Centr$(g)=\{h\in G\mid gh=hg\}$ denotes the centralizer of $g$. The following proposition, although elementary, will be very important to build actions of $G$ on $\Z^{k+1}$ in \S \ref{sec action by homeos}.

\begin{prop}\label{prop torsion free}
Let $G$ be a torsion-free and finitely-generated nilpotent group which is metabelian. Then:
\begin{itemize}
\item Given $g\in G$ and $0\neq m\in \mathbb{Z}$ we have that $\textup{Centr}(g^m)=\textup{Centr}(g).$
\item Let $A\leqslant G$ be a maximal abelian subgroup. If $ A$ is normal in $G$, then $G/A$ is torsion-free.
\end{itemize}
\end{prop}
\noindent\textbf{Proof.}
Assume that $f,g\in G$ and $m\in \mathbb{Z}$ are such that $[f,g^m]=e$. Define $H:=\langle f,g\rangle$, the subgroup generated by $f$ and $g$. Since $[g,g^m]=[f,g^m]=e$ we have that $g^m\in Z(H)$, and since $H/Z(H)$ is torsion-free (see  \cite[Proposition 5.2.19]{robinson}) we have that $g\in Z(H)$. Therefore, $\text{Centr}(g^m)\subseteq \text{Centr}(g)$ (the other inclusion is obvious).

The second point follows from the first. Let $A$ be a maximal abelian subgroup which is also normal, and assume $G/A$ is not torsion-free. Say $g\in G$ is such that $g\notin A$ but $g^m\in A$ for some $m\neq 0$. Then, since $A$ is abelian, we have that  $A\subseteq \text{Centr}(g^m)=\text{Centr}(g)$. In particular $\langle A, g\rangle$, the group generated by $A$ and $g$, is an abelian subgroup larger than $A$, contradicting our assumption. 

$\hfill \square$

\subsection{On the action of $G/A$ on $A$}\label{sec action G/A on A}
Let $G$ be a torsion-free and finitely-generated nilpotent group of degree $n$ which is also metabelian. Let $A$ be a maximal abelian subgroup containing $[G,G]$ (in particular it is normal). In view of Proposition \ref{prop torsion free}, we have that $G$ is an extension of $\mathbb{Z}^k$ by $\mathbb{Z}^d$,
$$1\longrightarrow \mathbb{Z}^d \longrightarrow G \longrightarrow \mathbb{Z}^k\longrightarrow 1,$$
\noindent where $A\simeq \mathbb{Z}^d$  and $G/A\simeq \mathbb{Z}^k$. In this section we study the natural action of $G/A$ on $A$ coming from the conjugacy action of $G$ on $A$.

 Let $\{ g_1,\ldots,g_d\}$ and  $\{ f_1A,\ldots,f_kA\}$ be generating sets of $A$ and $G/A$ respectively. Since $A$ is normal, the subgroup of $G$ generated by $f_1,\ldots,f_k$  acts on $A$ by automorphisms yielding a homomorphism $$\langle f_1,\ldots,f_k\rangle \longrightarrow \text{Aut}(\mathbb{Z}^d).$$ 
Therefore, the action of each $f\in \langle f_1,\ldots,f_k\rangle$ is given by a matrix ${A}_f\in GL_d(\mathbb{Z})$, which depends on the set $\{ g_1,\ldots,g_d\}$. We call $A_f$  the \emph{conjugacy matrix} of $f$. In the special case of the generators $f_1,\ldots, f_k$,  we denote the conjugacy matrix of $f_i$  simply by $A_i$.

In the next lemma, we will see that we can always choose a generating set of $A$, such that the conjugacy matrices of the elements $f_1,\ldots ,f_k$ belong to $U_d(\mathbb{Z})$, the group of upper triangular matrices with 1's in the diagonal. This is due to Mal'cev in the case where the matrix coefficients belong to a field. We write a direct proof in our special case. For the proof we will say that a generating set of a group is {\em minimal}, if it has least possible cardinality.

\begin{lemma}\label{lem triangular}
Let $A\leqslant G$ be a maximal abelian subgroup satisfying  that $[G,G]\subseteq A$. Suppose that $\mathbb{Z}^d\simeq A$ and $\mathbb{Z}^k\simeq G/A=\langle f_1A,\ldots,f_kA\rangle$.  
Then, there exists a generating set $\{g_1,\ldots,g_d\}$ of $A$ such that the conjugacy matrices of the elements $f_1,\ldots ,f_k$ belong to $U_d(\mathbb{Z})$.  In particular, the nilpotency degree of $G$ is bounded by $d+1$.
\end{lemma}

\noindent \textbf{Proof.}   Since $G$ is nilpotent of degree $n$,  the upper central series  
$$\{e\}=\zeta_0\leqslant \zeta_1\leqslant \cdots \leqslant \zeta_n=G,$$ 


\noindent is finite. Remember that all the factors $\zeta_i/\zeta_{i-1}$ are torsion-free. Combining this with the fact that $G/A$ is also torsion-free (see Proposition \ref{prop torsion free}), we have, for $g\in G$, that
\begin{equation}\label{eq torsion free}
g^j \in \zeta_i\cap A \Rightarrow g\in \zeta_i\cap A\hspace{0,3cm}\forall i \in \{0,\ldots,n\}, j\in \mathbb{Z}.
\end{equation}

\noindent Define $\Gamma_i:=\zeta_i\cap A$ and let $m$ be the smallest element in $\{1,...,n\}$ such that $\Gamma_m=A$. This yields the filtration $$\{e\}=\Gamma_0\leqslant \Gamma_1\leqslant \cdots \leqslant \Gamma_m=A,$$ 
such that 
 \begin{equation}\label {eq diagonal}[G,\Gamma_i]\subseteq \Gamma_{i-1},\end{equation}
 and, by (\ref{eq torsion free}),  also enjoys the property that each factor $\Gamma_{i+1}/\Gamma_i$ is torsion-free abelian. 
 
Note that if $\Gamma_{m-1}\simeq \Z^{n_{m-1}}$ and $\Gamma_m/\Gamma_{m-1}\simeq \Z^{n_m}$ then, since $\Gamma_m=A\simeq \Z^d$ is abelian, we have that $d=n_{m-1}+n_m$. Therefore, if $\{g_1,\ldots,g_{n_{m-1}}\}$ and $\{g_{n_{m-1}+1} \Gamma_{m-1},\ldots, g_{n_{m-1}+n_m}\Gamma_{m-1}\}$
 are minimal generating sets of $\Gamma_{m-1}$ and $\Gamma_m/\Gamma_{m-1}$ respectively, then 
 $\{g_1,\ldots,g_{n_{m-1}}, g_{n_{m-1}+1},\ldots , g_d\}$ is a minimal generating set of $\Gamma_m=A$. 
 
Recursively, we obtain a minimal generating set  $\{g_1,\ldots, g_d\}$ of $A$ which, by (\ref{eq diagonal}), has the property that for  $g_s\in \{g_1,\ldots,g_d\}\cap \Gamma_i$, it holds that $[f_j,g_s]\in \Gamma_{i-1}\subseteq \langle g_1,\ldots,g_{s-1}\rangle$ for all $ j \in \{1,\ldots,k\}.$
In other words the conjugacy matrices of each $f\in \{f_1,\ldots,f_k\}$ belong to $U_d(\Z)$. 

The fact that $G$ has nilpotency degree bounded by $d+1$ follows from the fact that $U_d(\Z)$ has nilpotency degree $d+1$.$\hfill\square$

\subsection{Invariant orders and their dynamical versions}
\label{sec orders}

We close these preliminaries with the concepts of order and dynamical realization. A group $G$ is left-orderable 
if it admits a total order relation, say $\preceq$, which is invariant under multiplication from the left, that is: $\text{if}\hspace{0,2cm}f\preceq g\hspace{0,2cm}\text{then}\hspace{0,2cm}hf\preceq hg\hspace{0,2cm}\text{for all}\hspace{0,2cm}h\in G.$
An important family of left-orderable groups are finitely-generated and torsion-free abelian groups. Indeed we will repetitively use the {\em lexicographic order} of $\Z^n$ defined by: 
\begin{equation} \label{eq lexicographic order} (i_1,\ldots,i_n)\prec (i_1^\prime,\ldots,i_n^\prime)\hspace{0,2cm} \Leftrightarrow \hspace{0,2cm}\exists\:k\in \{1,\ldots,n\}\hspace{0,2cm}\text{ such that }\hspace{0,2cm}i_k< i_k^\prime\hspace{0,2cm}\text{ and }\hspace{0,2cm}i_s=i_s^\prime\hspace{0,2cm}\text{ for }\hspace{0,2cm}s<k.\end{equation}

What is important for this work is that a countable group is left-orderable if and only if it embeds into  $\text{Diff}_+^0(\R)$ (see \cite[\S2]{na} or \cite[\S 1.1.3]{GOD} for details).  Since $\text{Diff}_+^0(\R)$ is isomorphic to $\text{Diff}_+^0([0,1])$, left-orderability of a countable group is equivalent to be isomorphic to a subgroup of $\text{Diff}_+^0([0,1])$. More generally, given a group $G$ that acts on a countable and totally ordered set $(\Omega,\preceq)$  by order preserving bijections, say $\omega \mapsto g(\omega)$,  for $g\in G$ and $\omega \in \Omega$, then there is a {\em dynamical realization} of this action. This means that there is an order preserving map $i:(\Omega,\preceq)\to ([0,1],\leq)$ and a homomorphism $\psi:G\to \text{Diff}_+^0([0,1])$ satisfying that
$\psi(g)(i(\omega))=i(g(\omega))$
for every $\omega\in \Omega$ and every $g\in G$. See \cite[Lemma 2.40]{BMRT} for a proof. Clearly $\psi$ is an embedding whenever the $G$ action on $\Omega$ is faithful.

\section{Proof of Theorem \ref{teo lower bound}}
\label{sec lower bound}

Throughout this section, $G$ will denote a non-abelian, torsion-free, finitely-generated nilpotent group which is metabelian, and $A$ will denote a maximal abelian subgroup containing $[G,G]$ (in particular it is normal). Recall that then $G$ is an extension of $G/A\simeq \Z^k$ by $A\simeq \Z^d$ (see \S \ref{sec action G/A on A}). In particular the nilpotency  degree of $G$ is bounded by $d+1$ (see Lemma \ref{lem triangular}).

\subsection{An action of $G$ on a totally ordered set}
\label{sec action by homeos}

\begin{prop}\label{prop polynomial bound} Fix a generating set $\{g_1,\ldots,g_d,f_1,\ldots,f_k\}$ of $G$, such that $\{ g_1,\ldots,g_d\}$ is a generating set of $A$ given by Lemma \ref{lem triangular} and  $\langle f_1A,\ldots,f_kA\rangle= G/A\simeq \mathbb{Z}^k$. Then, for a fixed $s\in\{1,\ldots,d\}$ there is an action of $G$ on $\mathbb{Z}^{k+1}$ satisfying:

\begin{enumerate} [1)]

\item   For all $m\in\{1,\ldots,d\}$ and all $t\in\{1,\ldots,k\}$ there exist functions $\ell_t,r_m:\mathbb{Z}^k\rightarrow \mathbb{Z}$, such that
$$f_t\cdot(i_1,..,i_t,..,i_k,j)=(i_1,..,i_t+1,..,i_k,j+\ell_t(i_1,\ldots,i_k)), $$ $$ g_m\cdot(i_1,\ldots,i_k,j)=(i_1,\ldots,i_k,j+r_m(i_1,\ldots,i_k)).$$ In particular,  the action of $G$ on $\mathbb{Z}^{k+1}$ preserves the lexicographic order.
 Besides, $r_s\equiv 1$ and $r_1=r_2=\cdots =r_{s-1}\equiv 0$.

\item There exists a positive constant $M$, such that for all $t\in \{1,\ldots,k\}$, $m\in \{1,\ldots,d\}$ and 

$(i_1,\ldots,i_k)\neq (0,\ldots,0)$ we have
$$\abs{\ell_t(i_1,\ldots,i_k)}\leq M(\abs{i_1}+\cdots+\abs{i_k})^d,\hspace{0,5cm} \abs{r_m(i_1,\ldots,i_k)}\leq M(\abs{i_1}+\cdots+\abs{i_k})^d.$$
\end{enumerate}
\end{prop}

\noindent \textbf{Proof.}  We start by showing item 1.  To this end, fix $s\in \{1,\ldots,d\}$ and consider the subgroup $H_s=\langle \{g_1,\ldots,g_d\}\smallsetminus \{g_s\} \rangle$. Since the sets $\{f_1^{i_1}\cdots f_k^{i_k} A:\;i_1,\ldots,i_k\in \mathbb{Z}\}$ and $\{g_s^jH_s:j\in \mathbb{Z}\}$ are partitions of $G$ and $A$ respectively, the coset space can be described by the {\em normal forms}
\begin{equation}\label{eq normal forms} G/H_s=\{f_1^{i_1}\cdots f_k^{i_k}g_s^jH_s\; : \;i_1,\ldots,i_k,j\in \mathbb{Z}\}.\end{equation}
Hence  we can identify $G/H_s$ with $\mathbb{Z}^{k+1}$ (as sets) by identifying  $f_1^{i_1}\cdots f_k^{i_k}g_s^jH_s$ with $({i_1},\ldots ,{i_k},j)$. In particular, the left-multiplication action of $G$ on $G/H_s$  provides an action of $G$ on $\Z^{k+1}$. This is the action we want to consider.

Now, by Lemma \ref{lem triangular}, we have that for all $i,j\in \{1,\ldots,k\}$ and $l\in \{1,\ldots,d\}$ it holds that  $$f_if_j\in f_jf_i\langle g_1,\ldots,g_d \rangle\hspace{0,2cm}\text{and}\hspace{0,2cm}g_lf_j \in f_jg_l\langle g_1,\ldots,g_{l-1}\rangle.$$ Therefore,  for $t\in \{1,\ldots, k\}$, the action of $f_t$ is addition by 1 on the $t$ coordinate and the action on the $k+1$ coordinate depends on previous $k$ coordinates, hence the function $\ell_t$. The function $r_m$, for $m\in \{1,\ldots,d\}$, can be found analogously. Finally, as the maps $\ell_t$ and $r_m$ depend only on the coordinates $(i_1,\ldots, i_k)$,  the reader can easily verify that the $G$ action on $\Z^{k+1}$ preserves the lexicographic order.

\vsp

Now we check item 2.  Let $t\in \{1,\ldots,k\}$. Recall that the action of $f_t$ on $\Z^{k+1}$ is nothing but the left-multiplication action of $f_t $ on $G/H_s$. Hence, in order to compute the image of $f_1^{i_1}\cdots f_t^{i_t}\cdots f_k^{i_k}g_s^jH_s$ under $f_t$, we need to multiply and  find the representative in normal form (\ref{eq normal forms}). 
To do this, observe that   $f_tf_j=[f_t,f_j]f_jf_t$. Hence, bringing $f_t$ to the $t$-th position generates at most 
$\abs{i_1}+\cdots+\abs{i_k}$ simple commutators of weight $2$, which we now need to move to the rightmost place ({\em i.e.} after the $f_k^{i_k}$ but before $g_s^j$).  Since $G$ is metabelian,  the commutators commute with each other. So, moving  them all to the rightmost place  generates at most $(\abs{i_1}+\cdots+\abs{i_k})^2$ simple commutators of weight $3$. Analogously, moving them all to the rightmost place, we have at most $(\abs{i_1}+\cdots+\abs{i_k})^3$ simple commutators of weight 4, and so on. Since $G$ has nilpotency degree bounded by $d+1$, all simple commutators of this weight are trivial (see Lemma \ref{lem triangular}). Therefore, repeating the previous argument $d+1$ times, we have 
$$f_t.(f^{i_1}_1\cdots f_t^{i_t}\cdots f_k^{i_k}g_s^jH_s)= f^{i_1}_1\cdots f_t^{i_t+1}\cdots f_k^{i_k}gg_s^jH_s,$$ 

\noindent where $g\in A$ is the product of at most $$\sum_{i=1}^d(\abs{i_1}+\cdots+\abs{i_k})^i \leq d(\abs{i_1}+\cdots+\abs{i_k})^d$$ simple commutators. Now note that $$g.g_s^jH_s=g_s^{\ell_t(i_1,\ldots,i_k)}g_s^jH_s,$$ since $\ell_t(i_1,\ldots,i_k)$ agrees with the exponent of $g_s$ in the expression of $g$ over the generators $g_1,\ldots,g_d$. Therefore, letting $\mathcal{S}\subseteq A$ be the set of all simple commutators of $G$ (which is finite), and defining $$\lambda:=\max\{\abs{m_s}: \exists\; m_1,\ldots,m_d\text{ for which }(g_1^{m_1}\cdots g_s^{m_s}\cdots g_d^{m_d})\in \mathcal{S}\},$$ 
we see that $\ell_t(i_1,\ldots,i_k)$ is bounded by $\lambda$ times the number of simple commutators that were used to write $g$. Hence $$\abs{\ell_t(i_1,\ldots,i_k)}\leq \lambda d(\abs{i_1}+\cdots+\abs{i_k})^d.$$ Analogous computations give the inequality for the functions $r_m$.

$\hfill\square$

\begin{remark}
Note that the action built in Proposition \ref{prop polynomial bound} is not necessarily faithful. However, it is such that the elements $g_1,\ldots,g_{s-1}$ act trivially and $g_s(i_1,\ldots,i_k,j)=(i_1,\ldots,i_k,j+1)$. This will be used in Section \ref{sec faithful action} in order to build a faithful action. 
\end{remark} 

\subsection{Action by diffeomorphisms of $[0,1]$}
\label{sec action by diffeos}
For a fixed  $s\in \{1,\ldots,d\}$, Proposition \ref{prop polynomial bound} builds an action  of $G$ on $\mathbb{Z}^{k+1}$ which preserves the lexicographic order. Hence we can consider the dynamical realization of this action (see the beginning of \S \ref{sec orders}) to get a $G$-action by orientation preserving homeomorphisms of $[0,1]$.
 
 Now, since the group is nilpotent  and we have good control from the polynomials appearing in  Proposition \ref{prop polynomial bound}, we will see that this action can actually be smoothed to an action by diffeomorphisms of $[0,1]$.  For this we need the following result from Pixton   and Tsuboi  \cite{pixton,tsuboi}. See the proof of Proposition 1.2 in \cite{tsuboi} for details.

\begin{lemma}\label{lem tsuboi} 

There exists a family of $C^{\infty}$-diffeomorphisms $\varphi_{I^\prime,I}^{J^\prime,J}:I\rightarrow J$, ranging over all bounded intervals $I$, $I'$, $J$, $J'$ of $\R$ where $I^\prime$ (resp. $J^\prime$) is adjacent to $I$ (resp. $J$) by the left, such that:

\begin{enumerate}[1)]

\item For all $I,I'$, $J,J',$  $K,K'$ as above,

$$ \varphi_{J^\prime,J}^{K^\prime,K}\circ\varphi_{I^\prime,I}^{J^\prime,J}=\varphi_{I^\prime,I}^{K^\prime,K}.$$

\item For all $I,I'$, $J,J'$, 
$$D\varphi^{J^\prime,J}_{I^\prime,I}(x_{-})=\frac{\abs{J^\prime}}{\abs{I^\prime}}\text{  and   }D\varphi^{J^\prime,J}_{I^\prime,I}(x_{+})=\frac{\abs{J}}{\abs{I}},$$ where $x_-$ (resp. $x_+$) is the left (resp. right) endpoint of $I$.

\item There is a constant $M$ such that for all  $I,I',J,J'$ as above, and all $x\in I$, we have
$$D\textup{log}(D\varphi_{I^\prime,I}^{J^\prime,J})(x)\leq \frac{M}{\abs{I}}\left\lvert\frac{\abs{I}\abs{J^\prime}}{\abs{J}\abs{I^\prime}}-1\right\rvert,$$ provided that $\max\{|I|, |I'|, |J|, |J'|\}\leq 2 \min\{|I|, |I'|, |J|, |J'|\}.$ 
\item Given $I, I', \;J,J', \;K,K',\: L,L'$,  as above, then 

$$\left\lvert \textup{log}(D\varphi_{I^\prime,I}^{K^\prime,K})(x)-\textup{log}(D\varphi_{J^\prime,J}^{L^\prime,L})(y)\right\rvert\leq \left\lvert\textup{log}\frac{\abs{K}\abs{J}}{\abs{I}\abs{L}}\right\rvert+ \left\lvert\textup{log}\frac{\abs{K^\prime}\abs{I}}{\abs{I^\prime}\abs{K}}\right\rvert + \left\lvert\textup{log}\frac{\abs{L^\prime}\abs{J}}{\abs{J^\prime}\abs{L}}\right\rvert,$$
for all $x\in I$, $y\in J$.
\end{enumerate}
\end{lemma}

\bigskip

Now, let $\{I_{i_1,\ldots,i_k,j}: (i_1,\ldots,i_k,j)\in \mathbb{Z}^{k+1}\}$ be a family of intervals  whose disjoint union is
dense in $[0, 1]$ and that are disposed preserving the lexicographic order of $\Z^{k+1}$. We identify the generators $g_1,\ldots,g_d,f_1,\ldots,f_k$ from  Lemma \ref{lem triangular} with elements in $\text{Diff}^0_+([0,1])$ as follows: 
$f_t$ and $g_s$ will be homeomorphisms of $[0,1]$ whose restriction to $I_{i_1,\ldots,i_k,j}$ coincides, respectively, with 

$$\varphi_{I_{i_1,...,i_k,j-1}\,,\, I_{i_1,...,i_k,j}}^{I_{i_1,..,i_t+1,..,i_k,j+\ell_t(i_1,...,i_k)-1}\,,\;I_{i_1,..,i_t+1,..,i_k,j+\ell_t(i_1,...,i_k)}}\hspace{0,3cm}\text{and}\hspace{0,3cm}\varphi_{I_{i_1,...,i_k,j-1},\,I_{i_1,...,i_k,j}}^{I_{i_1,...,i_k,j+r_s(i_1,...,i_k)-1}\,,\;I_{i_1,...,i_k,j+r_s(i_1,...,i_k)}},$$ for $t\in \{1,...,k\}$, and $s\in\{1,...,d\}$. 
Thus, by  1) in Lemma \ref{lem tsuboi},  we have a group homomorphism $G\rightarrow \text{Diff}^0_+([0,1]).$ The main technical step for proving Theorem A is the following proposition.

\begin{prop}\label{prop main} Given $\alpha<1/k$,  there is a choice of lengths of the intervals $|I_{i_1,\ldots,i_k,j}|$ such that  the homeomorphisms $f_1,\ldots,f_k,g_1,\ldots,g_d$ are  diffeomorphisms of class $C^{1+\alpha}$.

\end{prop}

The rest of Section \ref{sec action by diffeos} is devoted to the proof of Proposition \ref{prop main}. We will assume that $k\geq 2$ since, after Condition 3 in Proposition \ref{prop polynomial bound}, we can use the estimates from \cite[\S4]{int4} to ensure that, when $k=1$, the action is by $C^{1+\alpha}$ diffeomorphisms for any $\alpha <1$.

So let $k\geq 2$ and consider $\alpha<1/k$. Choose  positive real numbers  $p_1,\ldots,p_k,r$ such that for all $n\in \{1,\ldots,k\}$ the following conditions hold:
\begin{enumerate}[I)]
\item $\alpha + r\leq 2$,  
\item $d(r-1)\leq(1-\alpha)$,
\item $2dr\leq p_n$,
\item $2d\leq p_n(1-\alpha)$,
\item $\frac{1}{p_1}+\cdots+\frac{1}{p_k}+\frac{1}{r}< 1$,
\item $\alpha\leq \frac{1}{p_n}+\frac{1}{r}$ and $\alpha\leq \frac{r}{p_n(r-1)}$.
  
\end{enumerate} For instance, one can take $p_1=\cdots=p_k=3d/\alpha$ and $r=3d/(3d-1)$.

 Now define the lengths of the intervals $I_{i_1,\ldots,i_k,j}$ as 
$$\abs{I_{i_1,\ldots,i_k,j}}=\frac{1}{\abs{i_1}^{p_1}+\cdots+\abs{i_k}^{p_k}+\abs{j}^r+1}.$$  From condition V) it follows that $\sum \abs{I_{i_1,\ldots,i_k,j}}<\infty$, hence this family of intervals can be disposed on a finite interval respecting the lexicographic order. After renormalization, we can assume that this interval is $[0,1]$.

Following \cite{JNR}, we say that two real-valued functions $f$ and $g$ satisfy $f\prec g$ if there is a constant $M>0$ such that $\abs{f(x)}\leqslant Mg(x)$ for all $x$. We also write $f\asymp g$ if $f\prec g$ and $g\prec f$. 

 Let $\theta$ be a non-negative $C^2$ real-valued function satisfying $\theta(\xi)=\abs{\xi}^r$ for $\abs{\xi}\geqslant 1$, and $\theta(0)=0$. Consider the auxiliary functions ($C^2$ with respect to $\xi$):
 \begin{itemize}
 \item $\psi(i_1,\ldots,i_k,\xi):=1+\abs{i_1}^{p_1}+\cdots+\abs{i_k}^{p_k}+\theta(\xi),$  
 \item $\Psi_{i_1,\ldots,i_k}(\xi):=\text{log}(\psi(i_1,\ldots,i_k,\xi)).$
 \end{itemize}

 \begin{lemma} \label{lem asymptotic}
 Let $S=1+\abs{i_1}^{p_1}+\cdots+\abs{i_k}^{p_k}$. Given $C>0$, there exists a positive constant $M$ such that the inequality $$ \frac{1}{M}\psi(i_1,\ldots,i_k,j) \leq \psi(i_1,\ldots,i_k,\xi) \leq M \psi(i_1,\ldots,i_k,j),$$ holds for any $\xi$ satisfying $\abs{\xi-j}\leq C(S^{1/r}+(\abs{i_1}+\cdots +\abs{i_k})^d)$.

 \end{lemma} 
 
We remark that in the situation of Lemma \ref{lem asymptotic}, we will still use the notation $\psi(i_1,\ldots,i_k,j)\asymp\psi(i_1,\ldots,i_k,\xi)$.  Even if this is a slight abuse of notation, it is justified by comparing the functions $\psi_j$ and $\psi_\xi$ defined by $\psi_\xi (i_1,\ldots,i_k)=\psi (i_1,\ldots,i_k,\xi)$ and $\psi_j (i_1,\ldots,i_k)=\psi (i_1,\ldots,i_k,j)$, whenever the inequality $\abs{\xi-j}\leq C(S^{1/r}+(\abs{i_1}+\cdots +\abs{i_k})^d)$  holds.
\vsp\vsp

 \noindent\textbf{Proof of Lemma \ref{lem asymptotic}.} By symmetry, it is enough to show that 
$ \frac{\psi(i_1,\ldots,i_k,\xi)}{\psi(i_1,\ldots,i_k,j)}$  is bounded above. For this we note that
\begin{align*}
 \frac{\psi(i_1,\ldots,i_k,\xi)}{\psi(i_1,\ldots,i_k,j)}&\prec \frac{1+\abs{i_1}^{p_1}+\cdots+\abs{i_k}^{p_k}+\abs{j}^r+\abs{\xi-j}^r} {\psi(i_1,\ldots,i_k,j)}
\\& \prec 1+ \frac{S+(\abs{i_1}+\cdots+\abs{i_k})^{dr}}{\psi(i_1,\ldots,i_k,j)}\\&\prec 2+\frac{(\abs{i_1}+\cdots+\abs{i_k})^{dr}}{\psi(i_1,\ldots,i_k,j)},
\end{align*} 
where we repeatedly use the inequality $\abs{x+y}^a\prec \abs{x}^a+\abs{y}^a$, which holds for any $a>0$. Now just notice that the last expression is bounded. Indeed, since $(\abs{i_1}+\cdots+\abs{i_k})^{dr} \prec \abs{i_1}^{dr}+\cdots+\abs{i_k}^{dr}$, it is enough to observe that for each $n\in \{1,\ldots,k\}$, 
  $$\abs{i_n}^{dr}\leq (\psi(i_1,\ldots,i_k,j))^{\frac{dr}{p_n}}\leq\psi(i_1,\ldots,i_k,j),$$ 
which holds due to condition III). $\hfill\square$

\subsubsection{The maps $g_s$ are $C^{1+\alpha}$-diffeomorphisms }

We start the proof of Proposition \ref{prop main} by showing that the maps $g_s$, for $s\in \{1,\ldots,d\}$, are of class $C^{1+\alpha}$. That is, we want to show that $g_s$ is a $C^1$-diffeomorphism and that there is a constant $C>0$ such that 

\begin{equation}\label{eq Holder bound g}\notag \frac{\abs{\textup{log}Dg_s(x)-\textup{log}Dg_s(y)}}{\abs{x-y}^\alpha}\leq C\text{ for all different } x,y\in [0,1].
\end{equation}

To check this,  it is enough to find a uniform $C$ as above for  points $x,y$ in  $\cup_j I_{i_1,\ldots,i_k,j}$ (independent of $i_1,\ldots,i_k$). Indeed, after condition 2) in Lemma \ref{lem tsuboi} and the definition of $g_s$, it follows that $g_s$ has derivative 1 at the end points of the intervals $\cup_j I_{i_1,\ldots,i_k,j}$. Hence the conditions from \cite[Lemma 4.1.22]{na} are satisfied, therefore obtaining that the $g_s$'s are of class $C^{1+\alpha}$.

\vsp

\noindent {\bf Case 1:} The points $x,y$ belong to the same  $ I:=I_{i_1,\ldots,i_k,j}.$

Condition 3) in Lemma \ref{lem tsuboi} provides a Lipschitz constant for  $\log(Dg_s)$. So it is enough to bound
$$\frac{1}{\abs{I}^\alpha} \left| \frac{\abs{I}\abs{J^\prime}}{\abs{J}\abs{I^\prime}} -1\right|,$$  where $I^\prime=I_{i_1,\ldots,i_k,j-1}$, $J=I_{i_1,\ldots,i_k,j+r_s(i_1,\ldots,i_k)}$ and $J^\prime=I_{i_1,\ldots,i_k,j+r_s(i_1,\ldots,i_k)-1}$.

We will in fact bound the following (a posteriori) asymptotically equivalent expression $$\frac{1}{\abs{I}^\alpha}\text{log}\frac{\abs{I}\abs{J^\prime}}{\abs{J}\abs{I^\prime}}.$$
For this, notice that $\text{log}\frac{\abs{I}\abs{J^\prime}}{\abs{J}\abs{I^\prime}}$ is equal to 
$$\Psi_{i_1,\ldots,i_k}(j+r_s(i_1,\ldots,i_k))-\Psi_{i_1,\ldots,i_k}(j+r_s(i_1,\ldots,i_k)-1)-(\Psi_{i_1,\ldots,i_k}(j)-\Psi_{i_1,\ldots,i_k}(j-1)).
 $$ So applying the Mean Value Theorem first to the function $x\mapsto \Psi_{i_1,\ldots,i_k}(j+1+x)-\Psi_{i_1,\ldots,i_k}(j+x)$ and then to the function $x\mapsto D\Psi_{i_1,\ldots,i_k}(x)$ we have 
 \begin{equation}\label{eq to bound}  \left\lvert\text{log}\frac{\abs{I}\abs{J^\prime}}{\abs{J}\abs{I^\prime}}\right\rvert=\abs{r_s(i_1,\ldots,i_k)}\abs{D^2(\Psi_{i_1,\ldots,i_k})(\xi)}, 
  \end{equation} 
  where $\xi$ is a point in the convex hull of $\{j-1,j,j-1+r_s,j+r_s\}$.   Let us find an upper bound for $\abs{D^2(\Psi_{i_1,\ldots,i_k})(\xi)}$. Since $D\theta$ and $D^2\theta$ are bounded on $[-1,1]$, and
$$D^2(\Psi_{i_1,\ldots,i_k})(\xi)=\frac{D^2\theta(\xi)}{\psi(i_1,\ldots,i_k,\xi)}-\frac{(D\theta(\xi))^2}{(\psi(i_1,\ldots,i_k,\xi))^2},$$ 
 we have that  $$D^2(\Psi_{i_1,\ldots,i_k})(\xi)\prec \frac{1}{\psi(i_1,\ldots,i_k,\xi)}$$ for all $\xi\in [-1,1]$.  On the other hand for $\xi\notin [-1,1]$  we have that $\theta(\xi)=\abs{\xi}^r$. So, si nce   $\abs{\xi}^{r-2}<1$ and $\abs{\xi}^r/\psi(i_1,\ldots,i_k,\xi)<1$, it follows that
 \begin{equation}\label{eq cota derivada segunda} D^2(\Psi_{i_1,\ldots,i_k})(\xi)\prec \frac{\abs{\xi}^{r-2}}{\psi(i_1,\ldots,i_k,\xi)}\prec \frac{1}{\psi(i_1,\ldots,i_k,\xi)}. \end{equation}
Now going back to equation (\ref{eq to bound}) and using 3) of Proposition \ref{prop polynomial bound}, we have
$$\text{log}\frac{\abs{I}\abs{J^\prime}}{\abs{J}\abs{I^\prime}}\prec \frac{\abs{i_1}^d+\cdots+\abs{i_k}^d}{\psi(i_1,\ldots,i_k,\xi)}.$$ Note that for all $n\in \{1,\ldots,k\}$ condition IV) yields $$\abs{i_n}^d\leq (\psi(i_1,\ldots,i_k,\xi))^{\frac{d}{p_n}}\leq (\psi(i_1,\ldots,i_k,\xi))^{(1-\alpha)}.$$ Finally, thanks to the fact that $\xi$ belongs to the convex hull of $\{j-1,j,j-1+r_s,j+r_s\}$, we use the bounds of $r_s$ from  Proposition \ref{prop polynomial bound} to  apply  Lemma  \ref{lem asymptotic} and conclude that $$\frac{1}{\abs{I}^{\alpha}} \text{log}\frac{\abs{I}\abs{J^\prime}}{\abs{J}\abs{I^\prime}}\prec\frac{ (\psi(i_1,\ldots,i_k,\xi))^{-\alpha}}{\abs{I}^{\alpha}}\prec\frac{ (\psi(i_1,\ldots,i_k,j))^{-\alpha}}{\abs{I}^{\alpha}}=1,$$
as desired.

\noindent {\bf Case 2:} The point $x$ belongs to $I_{i_1,\ldots,i_k,j}$ and $y$ belongs to $I_{i_1,\ldots,i_k,j'}.$

\vsp

We assume without loss of generality that $j< j^\prime$. Condition 4) of Lemma \ref{lem tsuboi} tells us that $\abs{\textup{log}Dg_s(x)-\textup{log}Dg_s(y)}$ is bounded  above by 
$$\left\lvert\text{log}\frac{\abs{I_{i_1,\ldots,i_k,j+r_s}}\abs{I_{i_1,\ldots,i_k,j^\prime}}}{{\abs{I_{i_1,\ldots,i_k,j}}\abs{I_{i_1,\ldots,i_k,j^{\prime}+r_s}}}}\right\rvert+\left\lvert\text{log}\frac{\abs{I_{i_1,\ldots,i_k,j+r_s-1}}\abs{I_{i_1,\ldots,i_k,j}}}{{\abs{I_{i_1,\ldots,i_k,j-1}}\abs{I_{i_1,\ldots,i_k,j+r_s}}}}\right\rvert+\left\lvert\text{log}\frac{\abs{I_{i_1,\ldots,i_k,j^\prime+r_s-1}}\abs{I_{i_1,\ldots,i_k,j^\prime}}}{{\abs{I_{i_1,\ldots,i_k,j^\prime-1}}\abs{I_{i_1,\ldots,i_k,j^\prime+r_s}}}}\right\rvert.$$ 

The estimates in Case 1 allow us to control the last two terms (divided by $ |x-y|^\alpha$), thus we only need to bound the first term. So we look for a uniform bound for 
\begin{equation}\label{eq to bound 2}  \frac{1}{\abs{x-y}^\alpha}\left\lvert\text{log}\frac{\abs{I}\abs{J^\prime}}{{\abs{I^\prime}\abs{J}}}\right\rvert,\end{equation} where $I=I_{i_1,\ldots,i_k,j}$, $I^\prime=I_{i_1,\ldots,i_k,j^\prime}$, $J=I_{i_1,\ldots,i_k,j+r_s}$ and  $J^\prime=I_{i_1,\ldots,i_k,j^\prime+r_s}$. Assume that $j$, $j^\prime$ are positive (the case where both are negative follows by symmetry, and if they have different
sign, it suffices to consider an intermediate comparison with the term corresponding to $j^{\prime\prime}=0$). Assume further that $j^\prime-j\geq 2$ (the case where $j^\prime-j= 1$ follows from Case 1, passing through the point that separates the intervals and using triangle inequality). 
Again, applying the Mean Value Theorem first to the function $x\mapsto \Psi_{i_1,\ldots,i_k}(j+1+x)-\Psi_{i_1,\ldots,i_k}(j+x)$  and then to the function $x \mapsto D\Psi_{i_1,\ldots,i_k}(x)$ we have
\begin{equation} \label{eq bound derivative}\left\lvert\text{log}\frac{\abs{I}\abs{J^\prime}}{{\abs{I^\prime}\abs{J}}}\right\rvert=\abs{j-j^\prime}\cdot\abs{r_s(i_1,\ldots,i_k)}\cdot\abs{D^2(\Psi_{i_1,\ldots,i_k})(\xi)},\end{equation} for a certain $\xi$ in the convex hull of $\{j,j^\prime,j+r_s,j^\prime+r_s\}$.

We start by bounding $|x-y|^{-\alpha}$. For this note that by Case 1 and the triangle inequality, we can (and will)  assume that $x$ is the left endpoint of $I$ and $y$ is the right endpoint of $I^\prime$. This yields 
$$\frac{1}{\abs{x-y}^\alpha}=\left(\frac{1}{\sum_{\ell=j}^{j^\prime}\abs{I_{i_1,\ldots,i_k,\ell}}}\right)^\alpha\leq \left(\frac{1}{\abs{j-j^\prime}\abs{I_{i_1,\ldots,i_k,j^\prime}}}\right)^\alpha,$$ where the last inequality holds because $\abs{I_{i_1,\ldots,i_k,j^\prime}}<\abs{I_{i_1,\ldots,i_k,\ell}}$ for $\ell<j^\prime$. Note that if in addition $\abs{j^\prime-j}\leq C( S^{1/r}+(\abs{i_1}+\cdots+\abs{i_k})^d)$, for some $C>0$, we can use Lemma \ref{lem  asymptotic} to compare $\abs{I}$ with $\abs{I^\prime}$, and eventually obtain the inequality

\begin{equation}\label{eq 5}
\frac{1}{\abs{x-y}^\alpha}\prec \left(\frac{1}{\abs{j-j^\prime}\abs{I_{i_1,\ldots,i_k,j}}}\right)^\alpha.
\end{equation}

We now exhibit a bound for (\ref{eq to bound 2}). We consider three  separate cases. Let $M$ be the constant in Proposition \ref{prop polynomial bound}.

\begin{enumerate}
\item[$i)$] The integers $j,j^\prime$ belong to $[0,2M(\abs{i_1}+\cdots+\abs{i_k})^d]$. Since $\xi\in \text{conv}\{j,j^\prime, j+r_s,j^\prime+r_s\}$   it follows from (\ref{eq cota derivada segunda}) and Lemma \ref{lem asymptotic} that 
$$\abs{D^2(\Psi_{i_1,\ldots,i_k})(\xi)}\prec\frac{1}{\psi(i_1,\ldots,i_k,\xi)}\asymp \frac{1}{\psi(i_1,\ldots,i_k,j)}.$$ Furthermore, we have that
$$\abs{j-j^\prime}\abs{r_s(i_1,\ldots,i_k)} \prec (\abs{i_1}+\cdots+\abs{i_k})^{2d}\prec (\psi(i_1,\ldots,i_k,j))^{1-\alpha},$$ where the last inequality holds by condition IV). If we combine this with  (\ref{eq bound derivative}), (\ref{eq 5}),  we conclude that
 $$\frac{1}{\abs{x-y}^\alpha}\left\lvert\text{log}\frac{\abs{I}\abs{J^\prime}}{{\abs{I^\prime}\abs{J}}}\right\rvert \prec \frac{1}{\abs{I}^\alpha}\frac{(\psi(i_1,\ldots,i_k,j))^{1-\alpha}}{\psi(i_1,\ldots,i_k,j)}= 1.$$
 
\item[$ii)$] The integers $j,j'$ belong to\footnote{the constant $k^d$ is just to ensure that the interval is non-empty. }  $ [2M(\abs{i_1}+\cdots+\abs{i_k})^d, 2M k^d S^{1/r}]$. Similarly to $i)$, the reader can check that we are in the hypotheses of Lemma \ref{lem asymptotic} and that $\abs{\xi}\geq M(\abs{i_1}+\cdots+\abs{i_k})^d$. Therefore, by  (\ref{eq cota derivada segunda}), (\ref{eq bound derivative}) and (\ref{eq 5}), we get
\begin{align*}
\frac{1}{\abs{x-y}^\alpha}\left\lvert\text{log}\frac{\abs{I}\abs{J^\prime}}{{\abs{I^\prime}\abs{J}}}\right\rvert &\prec \left(\frac{1}{\abs{j-j^\prime}\abs{I_{i_1,\ldots,i_k,j}}}\right)^\alpha \abs{j^\prime-j}\frac{(\abs{i_1}+\cdots+\abs{i_k})^d\abs{\xi}^{r-2}}{\psi(i_1,\ldots,i_k,\xi)} \\& \prec
\abs{j^\prime-j}^{1-\alpha}\frac{(\abs{i_1}+\cdots+\abs{i_k})^{d(r-1)}}{\psi(i_1,\ldots,i_k,j)^{1-\alpha}}.
 \end{align*}
To prove that this last expression is bounded, it is enough to show that 
 $$\abs{j^\prime-j}^{1-\alpha}(\abs{i_1}+\cdots+\abs{i_k})^{d(r-1)}\prec \psi(i_1,\ldots,i_k,j)^{1-\alpha}.$$ Since $j^\prime-j\leq 2M k^d S^{1/r}$, it follows that
 $$\abs{j^\prime-j}^{1-\alpha}(\abs{i_1}+\cdots+\abs{i_k})^{d(r-1)}\prec (1+\abs{i_1}^{p_1}+\cdots+\abs{i_k}^{p_k})^{\frac{(1-\alpha)}{r}}(\abs{i_1}+\cdots+\abs{i_k})^{d(r-1)},$$
so it suffices to prove that, given $n,m\in \{1,...,k\}$ we have
\begin{equation}\label{eq final ii}
\abs{i_n}^{\frac{p_n(1-\alpha)}{r}}\abs{i_m}^{d(r-1)}\prec (\psi(i_1,...,i_k,j))^{1-\alpha}.
\end{equation}
But  note that  $$\abs{i_n}^{\frac{p_n(1-\alpha)}{r}}\abs{i_m}^{d(r-1)}\leq (\psi(i_1,\ldots,i_k,j))^{\frac{(1-\alpha)}{r}+\frac{d(r-1)}{p_m}},$$ 
and that conditions II) and V) guarantee $\frac{(1-\alpha)}{r}+\frac{d(r-1)}{p_m}\leq (1-\alpha)$, which implies (\ref{eq final ii}).

\item[$iii)$] Finally suppose that the integers $j,j^\prime$ belong to $ [2M k^d S^{1/r},\infty)$.

 If $j^\prime\leq 2j$, then 
\begin{equation}
\frac{\psi(i_1,\ldots,i_k,j^\prime)}{\psi(i_1,\ldots,i_k,j)}\prec 1+\frac{\abs{j-j^\prime}^r}{\psi(i_1,\ldots,i_k,j)}\prec 1+ \frac{\abs{j}^r}{\psi(i_1,\ldots,i_k,j)}\leq 2. 
\end{equation} 
In particular, the intervals  $\abs{I^\prime}$  and $\abs{I}$ have comparable size and hence we conclude that (\ref{eq 5}) still holds. Also note that $j^\prime\leq 2j$ implies $\abs{\xi-j}\leq \abs{j}+M(\abs{i_1}+\cdots+\abs{i_k})^d$. Then, proceeding as in $ii)$, we have that
$$\frac{1}{\abs{x-y}^\alpha}\left\lvert\text{log}\frac{\abs{I}\abs{J^\prime}}{\abs{I^\prime}\abs{J}}\right\rvert\prec \frac{\abs{j}^{1-\alpha}(\abs{i_1}+\cdots+\abs{i_k})^{d(r-1)}}{\psi(i_1,\ldots,i_k,j)^{1-\alpha}}.$$
The reader can check, again as in $ii)$, that this last expression is bounded.

\bigskip

For the case $j^\prime> 2j$. We have 
\begin{align*}
\abs{x-y}&=\sum_{\ell=j}^{j^\prime}\abs{I_{i_1,\ldots,i_k,\ell}}=\sum_{\ell=j}^{j^\prime}\frac{1}{\abs{i_1}^{p_1}+\cdots+\abs{i_k}^{p_k}+\abs{\ell}^r}\\& \succ \sum_{\ell=j}^{j^\prime} \frac{1}{\abs{\ell}^r}\succ \int_{\ell=j}^{j^\prime} \frac{1}{x^r} dx\succ \frac{1}{\abs{j}^{r-1}},
\end{align*}
where the last inequality holds because $j'>2j$. On the other hand, applying the Mean Value Theorem, it follows that
\begin{align*}
\text{log}\frac{\abs{I}\abs{J^\prime}}{\abs{I^\prime}\abs{J}}&=\abs{r_s(i_1,\ldots,i_k)}\abs{D(\Psi_{i_1,\ldots,i_k})(\xi)-D(\Psi_{i_1,\ldots,i_k})(\tilde{\xi})},
\end{align*}
where $\xi\in \text{conv}\{j,j+r_s\}$ and $\tilde{\xi}\in \text{conv}\{j^\prime,j^\prime+r_s\}$. Therefore, observing that the function $\xi\mapsto D(\Psi_{i_1,\ldots,i_k})(\xi)=r\abs{\xi}^{r-1}/\psi(i_1,\ldots,i_k,\xi)$ is decreasing, we have
$$\frac{1}{\abs{x-y}^\alpha}\left\lvert\text{log}\frac{\abs{I}\abs{J^\prime}}{\abs{I^\prime}\abs{J}}\right\rvert\prec (\abs{i_1}+\cdots+\abs{i_k})^d\frac{\abs{j}^{(\alpha+1)(r-1)}}{\psi(i_1,\ldots,i_k,j)}.$$
Now we want to see that this last expression is bounded, in other words that the inequality $\abs{j}^{(\alpha +1)(r-1)}(\abs{i_1}+\cdots+\abs{i_k})^d\prec \psi(i_1,\ldots,i_k,j)$ holds. For this, arguing as in (\ref{eq final ii}), it is enough to check that for all  $ n\in \{1,...,k\}$ the inequality  
 $$\frac{(\alpha +1)(r-1)}{r}+\frac{d}{p_n}\leq 1$$ holds.
To see this, note that from IV) it follows that $\frac{d}{p_n}\leq \frac{1-\alpha}{2}$. Finally notice that $\frac{(\alpha +1)(r-1)}{r}+\frac{1-\alpha}{2}\leq 1\Leftrightarrow r \leq 2$, which is ensured by condition I).

\end{enumerate}

\subsubsection{ The maps $f_t$ are $C^{1+\alpha}$-diffeomorphisms  }

In the same way that for the maps $g_s$, we want to see that there is a constant $C > 0$ such that
$$\frac{\abs{\text{log}Df_t(x)-\text{log}Df_t(y)}}{\abs{x-y}^\alpha}\leqslant C \text{ for all different } x,y\in [0,1].$$
To simplify notation we will  only work with $t=1$ as the other cases are analogous.  As for the case of the maps $g_s$, we only have two cases to analyse.

\vsp

\noindent {\bf Case 1:} The points $x,y$ belongs to the same interval $I_{i_1,\ldots,i_k,j}$.

By Lemma \ref{lem tsuboi} it is enough to show that the following expression is uniformly bounded $$\frac{1}{\abs{I_{i_1,\ldots,i_k,j}}^\alpha}\text{log}\frac{\abs{I_{i_1,\ldots,i_k,j}}\abs{I_{i_1+1,i_2,\ldots,i_k,j+\ell_1-1}}}{\abs{I_{i_1+1,i_2,\ldots,i_k,j+\ell_1}}\abs{I_{i_1,\ldots,i_k,j-1}}}.$$ To see this, simply note that the above expression is equal to 
$$\frac{1}{\abs{I_{i_1,\ldots,i_k,j}}^\alpha}\text{log}\frac{\abs{I_{i_1,\ldots,i_k,j}}\abs{I_{i_1+1,i_2,\ldots,i_k,j-1}}}{\abs{I_{i_1+1,i_2,\ldots,i_k,j}}\abs{I_{i_1,\ldots,i_k,j-1}}}+
\frac{1}{\abs{I_{i_1,\ldots,i_k,j}}^\alpha}\text{log}\frac{\abs{I_{i_1+1,i_2,\ldots,i_k,j}}\abs{I_{i_1+1,i_2,\ldots,i_k,j+\ell_1-1}}}{\abs{I_{i_1+1,i_2,\ldots,i_k,j+\ell_1}}\abs{I_{i_1+1,i_2,\ldots,i_k,j-1}}}.$$
By condition VI)  we know from  \cite[\S 3.3]{int4} that the first term is uniformly bounded.  The second term is bounded as well since it is the same that we bounded when dealing with $g_s$ (changing $i_1$ for $i_1+1$).

\vsp

\noindent {\bf Case 2:} The point  $x\in I=I_{i_1,\ldots,i_k,j}$ and $y\in J=I_{i_1,\ldots,i_k,j^\prime}$, with $j< j^\prime$.

Here  we can use 4) from Lemma \ref{lem tsuboi} to bound $\abs{ \log Df_1(x)-\log Df_1(y)}$ by 
$$\left\lvert\text{log}\frac{\abs{I_{i_1+1,\ldots,i_k,j+\ell_1}}\abs{I_{i_1,\ldots,i_k,j^\prime}}}{{\abs{I_{i_1,\ldots,i_k,j}}\abs{I_{i_1+1,\ldots,i_k,j^{\prime}+\ell_1}}}}\right\rvert+\left\lvert\text{log}\frac{\abs{I_{i_1+1,\ldots,i_k,j+\ell_1-1}}\abs{I_{i_1,\ldots,i_k,j}}}{{\abs{I_{i_1,\ldots,i_k,j-1}}\abs{I_{i_1+1,\ldots,i_k,j+\ell_1}}}}\right\rvert+\left\lvert\text{log}\frac{\abs{I_{i_1+1,\ldots,i_k,j^\prime+\ell_1-1}}\abs{I_{i_1,\ldots,i_k,j^\prime}}}{{\abs{I_{i_1,\ldots,i_k,j^\prime-1}}\abs{I_{i_1+1,\ldots,i_k,j^\prime+\ell_1}}}}\right\rvert, $$  
and then work in the same way as for the functions $g_s$. For example, we express the term
 $$\frac{1}{\abs{x-y}^\alpha}\text{log}\frac{\abs{I_{i_1+1,i_2,\ldots,i_k,j+\ell_1}}\abs{I_{i_1,\ldots,i_k,j^\prime}}}{{\abs{I_{i_1,\ldots,i_k,j}}\abs{I_{i_1+1,i_2,\ldots,i_k,j^{\prime}+\ell_1}}}}$$  as
$$\frac{1}{\abs{x-y}^\alpha}\text{log}\frac{\abs{I_{i_1+1,i_2,\ldots,i_k,j}}\abs{I_{i_1,\ldots,i_k,j^\prime}}}{{\abs{I_{i_1+1,i_2,\ldots,i_k,j^\prime}}\abs{I_{i_1,\ldots,i_k,j}}}}+\frac{1}{\abs{x-y}^\alpha}\text{log}\frac{\abs{I_{i_1+1,i_2,\ldots,i_k,j+\ell_1}}\abs{I_{i_1+1,i_2,\ldots,i_k,j^\prime}}}{{\abs{I_{i_1+1,i_2,\ldots,i_k,j}}\abs{I_{i_1+1,i_2,\ldots,i_k,j^{\prime}+\ell_1}}}}.$$  The first term is bounded  by \cite[\S 3.3]{int4} and the second is also bounded by the same argument used for the functions $g_s$.

\subsection{Faithful actions}
\label{sec faithful action}

Given $s\in\{1,\ldots,d\}$ and a compact interval $I_s$, we have seen in Proposition \ref{prop main} how to produce an  action
 $$\phi_s : G \rightarrow \text{Diff}_+^{1+\alpha}(I_s).$$ 
Recall that the action from Proposition \ref{prop main} is a smoothing of the dynamical realization of the action given in Proposition  \ref{prop polynomial bound}. In particular,  the subgroup $\langle g_1,\ldots,g_{s-1}\rangle$ acts trivially, while $\langle g_s\rangle$ acts faithfully.
 
To obtain a faithful action of $G$  we do the following:
consider compact intervals $I_1,\ldots,I_d$ such that for all $s\in \{1,\ldots,d-1\}$, $I_{s+1}$ is contiguous to $I_s$ by the right. Then define on $I:=I_1\cup \cdots \cup I_d$ the action $\phi: G \rightarrow \text{Diff}_+^{1+\alpha}(I)$ as $$\phi\mid_{I_s}=\phi_s.$$ We claim that $\phi$ is injective. Indeed, since $Z(G)\leqslant A=\langle g_1,\ldots, g_d\rangle$, by  Proposition \ref{prop inyectiva} we only need to check that $\phi\mid_A$ is injective. Let $g\in A$ be an element that acts trivially on $I$. Then,  there exist $j_1,\ldots,j_d\in \mathbb{Z}$ such that $g=g_1^{j_1}\cdots g_d^{j_d}$.
Now, since $\phi(g)=id$, it follows that
$$\phi_s(g)=id,\; \forall s\in \{1,\ldots,d\}.$$ This yields that $j_d=\cdots=j_1=0$ and hence $g$ is the trivial element. This finishes the proof of Theorem A.

\section{Examples}
\label{sec examples}

In this section we give  examples of nilpotent groups for which we can compute the critical regularity.  In each case we use Theorem \ref{teo lower bound} to obtain a lower bound for the critical regularity and we argue that in our examples this is also an upper bound for the regularity.

We begin by recalling that if $G$ is a finitely-generated nilpotent group of homeomorphisms of $(0,1)$ that has no global fixed points, then there is a well-defined group homomorphism $\rho :G \rightarrow \mathbb{R},$ which is usually called the {\em translation number} of the action.  This map characterizes the elements of $G$ that have fixed points, in the sense that $\rho(g)=0$ if and only if $g$ has a fixed point in $(0,1)$.   Further, the action of $G$ on the interval has {\em no crossings}. By this we mean  that if an element $f\in G$ fixes an open subinterval  $I$ of $(0,1)$  and satisfies that $f(x)\neq x \text{ for all } x \text{ in }I$,  then for  any other $g\in G$ we have that $g(I)=I$ or $g(I)\cap I=\emptyset$. See \cite[\S 2.2.5]{na} for details. With this, it is easy to prove the following result that we will repeatedly use. 
\begin{lemma}\label{lem conrad} 
Let $G\leqslant \textup{Diff}_+^0(0,1)$ be a nilpotent group, and let $c\in G$ be a non-trivial element such that $c=[a,b]$  for some elements $a,b\in G$. If $c$ fixes an open interval $I$ and has no fixed point inside, then either $a$ or $b$ moves $I$ (disjointly). 
\end{lemma}
\noindent\textbf{Proof.} Looking for a contradiction, assume that $a$ and $b$ fix $I$. Then we have the translation number homomorphism for the group $\langle a, b, c\rangle\leqslant \text{Diff}_+^0(I)$. Since $c$ is a commutator, it is in the kernel of this homomorphism. Hence we conclude that $c$ has a fixed point inside $I$, which is contrary to our assumptions.
$\hfill\square$

To obtain upper bounds for the regularity of our groups we will  use a  result from Deroin, Kleptsyn and Navas \cite{DKN_acta}. We use the version from \cite[Proposition 2.1]{int4}.

\begin{theorem} \label{teo control of distortion}
Let $f_1,\ldots,f_k$ be $C^1$-diffeomorphisms of the interval $[0, 1]$ that commute with a $C^1$-diffeomorphism $g$. Assume that $g$ fixes a subinterval $I$ of $[0, 1]$ and its restriction to $I$ is non-trivial. Assume moreover that for a certain $0 < \alpha < 1$ and a sequence of indexes $i_j\in \{1,\ldots,k\}$, the sum
$$\sum_{j\geq 0}\abs{f_{i_j}\cdots f_{i_1}(I)}^\alpha <\infty. $$ Then $f_1,\ldots, f_k $ cannot be all of class $C^{1+\alpha}$.  

\end{theorem}
The following lemma is useful to get into the hypotheses of Theorem \ref{teo control of distortion}. Although it is stated in a slightly different way, the reader can check that the proof is exactly the same as that of \cite[Lemma 3.3]{DKN_acta}. 

\begin{lemma}\label{lem caminatas} 
Let $f_1,\ldots,f_k$ be $C^1$-diffeomorphisms of $[0,1]$, and $I$ a subinterval of $[0,1]$ such that $\mathbb{Z}^k\simeq \langle f_1,\ldots ,f_k\rangle/\textup{Stab}(I)$, where $\textup{Stab}(I)$ is the stabilizer of $I$ (which is assumed to be a normal subgroup). Then, if $\alpha> 1/k$ there exists a sequence $(f_{i_j})_{j\in \mathbb{N}}$ of elements in $\{f_1,\ldots ,f_k\}$ such that $$\sum_{j\geq 0}\abs{f_{i_j}\cdots f_{i_1}(I)}^\alpha< \infty .$$
\end{lemma}

\subsection{Heisenberg Groups}
\label{sec Heissemberg}


For a natural number $n\geq 1$, the discrete $(2n+1)$-dimensional \emph{Heisenberg group}, is defined as the set of matrices 
$$\mathscr{H}_n:=\left\lbrace\begin{pmatrix}
1 & \vec{x} & c\\
\vec{0}^t & I_n & \vec{y}^{\, t}\\
0 & \vec{0} & 1
\end{pmatrix}: \vec{x},\vec{y}\in \mathbb{Z}^n, c\in \mathbb{Z} \text{ and } I_n\text{ is the identity matrix of size }n\right\rbrace,$$
with the usual matrix product.  Note that the center of $\mathscr H_n$ coincides with the commutator subgroup and is generated by the matrix
$$\textbf{C}:=
\begin{pmatrix}
1 & \vec{0} & 1\\
\vec{0}^t & I_n & \vec{0}^{\, t}\\
0 & \vec{0} & 1
\end{pmatrix}.$$

We want to prove Theorem \ref{teo Heissemberg}, but before this, it will be useful for us to bound the rank of  maximal abelian subgroups. Assume that there exists a maximal abelian subgroup of $ \mathscr{H}_n$ of rank $m$. Then we can choose elements 
$$\textbf{A}_i:=\begin{pmatrix}
1 & \vec{a_i} & c_i\\
\vec{0}^t & I_n & \vec{b_i}^t\\
0 & \vec{0} & 1
\end{pmatrix}\in \mathscr{H}_n\hspace{0,2cm}\text{for}\hspace{0,2cm}i\in \{1,\ldots,m-1\},$$
such that $\langle \textbf{A}_1,\ldots \textbf{A}_{m-1},\textbf{C}\rangle\simeq \mathbb{Z}^{m}.$ Note that the commutativity of these matrices is equivalent to the equations
\begin{equation}\label{eq matrix}
\vec{a_i}\cdot\vec{b_j}=\vec{a_j}\cdot\vec{b_i}\hspace{0,3cm}\forall\; i,j\in \{1,...,m-1\}.
\end{equation} 
Note also that if $\{\textbf{A}_1,\ldots, \textbf{A}_{m-1},\textbf{C}\}$ generates a free abelian group of rank $m$, then the set of vectors $\mathscr{B}:=\{(\vec{b}_i,\vec{a}_i)\in \mathbb{Z}^{n}\times \mathbb{Z}^{n}:1\leq i\leq m-1\}$ generates a free abelian group of rank $m-1$. Indeed, if we have a dependency relation, say 
$r(\vec{b}_1,\vec{a}_1)\in \langle (\vec{b}_2,\vec{a}_2),\ldots,(\vec{b}_{m-1},\vec{a}_{m-1})\rangle$ for some $0\neq r\in\mathbb{Z},$ then $\textbf{A}_1^{r}\in \langle \textbf{A}_2,\ldots,\textbf{A}_{m-1},\textbf{C}\rangle,$ which contradicts that the abelian group has rank $m$. 

Having said this, we claim that $m\leq n+1$. To see the latter, note that, by equations (\ref{eq matrix}) any vector of the form $(\vec a_i,-\vec b_i)$, with $1\leq i\leq m -1$, is perpendicular to $\langle \mathscr B\rangle$. Hence we have two orthogonal subgroups of $\Z^n\times \Z^n$ that both have rank $m-1$, and thus $m-1\leq n$, which proves our claim.

\paragraph{Realization.}
Consider the abelian subgroup 
$$A:=\left\lbrace \begin{pmatrix}
1 & \vec{x} & c\\
\vec{0}^{\,t} & I_n & \vec{0}^{\,t}\\
0 & \vec{0} & 1
\end{pmatrix}\;: \vec{x}\in \mathbb{Z}^n\text{ and }c\in\mathbb{Z} \right\rbrace.$$ 
Notice that $A$ has rank equal to $n+1$, which is the largest we can expect. Since the rank of $\mathscr{H}_n/A$ is $n$, we have that Theorem \ref{teo lower bound} provides an injective group homomorphism
$$\mathscr{H}_n\hookrightarrow \text{Diff}_+^{1+\alpha}([0,1])\text{ for } \alpha < 1/n.$$
\paragraph{Bounding the regularity.}
Now we consider a faithful action $\phi:\mathscr{H}_n\hookrightarrow \text{Diff}_+^1([0,1])$. Abusing notation, we can think $\mathscr{H}_n\leqslant \text{Diff}_+^1([0,1])$. 

Since  the commutator subgroup of $\mathscr{H}_n$ is generated by ${\bf C}$, we deduce from Lemma \ref{lem conrad} that {\bf C} has fixed points inside $(0,1)$. Therefore, we can find an interval
$I\subsetneq [0,1]$ such that
$\textbf{C}(I)=I\text{ and }\textbf{C}(x)\neq x$ for all $x$ in the interior of  $I$. Let $\text{Stab}(I)$ be the stabilizer of $I$. It is easy to see that this is an abelian subgroup. Indeed, if we take $\textbf{A},\textbf{B}\in \text{Stab}(I)$ and assume that they do not commute, then there must exist $m\in \mathbb{Z}$ such that $[\textbf{A},\textbf{B}]=\textbf{C}^m$. Since $\textbf{C}$ has no fixed points inside $I$, Lemma \ref{lem conrad} tell us that either $\textbf{A}$ or $\textbf{B}$ moves $I$ , which is a contradiction.  Note that Stab($I$) is a normal subgroup as it contains the commutator subgroup. Further, we know that there is a natural number $k$ and elements $\textbf{B}_1,\ldots,\textbf{B}_k\in \mathscr{H}_n$ such that  $$\mathbb{Z}^k\simeq\frac{\mathscr{H}_n}{\text{Stab}(I)}= \frac{\langle\textbf{B}_1,\ldots,\textbf{B}_k\rangle}{\text{Stab}(I)}.$$ So, given $\alpha> 1/k$,  we can find by Lemma \ref{lem caminatas} a sequence $(\textbf{B}_{i_j})_{j\in\mathbb{N}}$ of elements in $\{\textbf{B}_1,\ldots,\textbf{B}_k\}$ such that 
$$\sum_{j\geq 0}\abs{\textbf{B}_{i_j}\cdots \textbf{B}_{i_1}(I)}^\alpha<\infty,$$
and hence  Theorem \ref{teo control of distortion} yields that $\phi$ is not an action by $C^{1+\alpha}$-diffeomorphisms. 

Now since the rank of $\text{Stab}(I)$ is bounded above by $n-1=\text{rank}(A)$, we have that  $$k=\text{rank}\left(\frac{\mathscr{H}_n}{\text{Stab}(I)}\right)\geq \text{rank}\left(\frac{\mathscr{H}_n}{A}\right)=n, $$
which implies that the regularity of the action $\phi$ is bounded above by $1+1/n$.  So we conclude that
$$\text{Crit}_{[0,1]}(\mathscr{H}_n)=1+\frac{1}{n}.$$

\begin{remark}
If $G$ is a finitely-generated, torsion-free nilpotent group whose center is cyclic and satisfies $[G,G]\leqslant Z(G)$, then the proof of Theorem \ref{teo Heissemberg} yields  that the lower bound for $\textup{Crit}_{[0,1]}(G)$ given by  Theorem \ref{teo lower bound} is the critical one. 

\end{remark}

\subsection{Examples with large nilpotency degree}

\label{sec more examples} 
Theorem \ref{teo Heissemberg} gives us the critical regularity for the Heisenberg groups, which are groups having nilpotency degree 2. In this section we provide more examples of  nilpotent groups for which we can compute the critical regularity, but whose nilpotency degree can be arbitrarily large.  As for the Heisenberg groups, in these examples we will show that the lower bound provided by Theorem  \ref{teo lower bound} is also an upper bound.

 Fix $d,k \in \mathbb{N}$, assume $d\geq k$ and consider a matrix $(m_{i,s})\in M_{k}(\mathbb{Z})$ with non-zero determinant and  positive entries. We let $G$ be the group generated by the set
$$\{g_0\}\cup \{g_{i,j}: (i,j)\in \{1,\ldots,k\}\times\{1,\ldots,d\}\}\cup \{f_1,\ldots,f_k\} ,$$  subject to the relations

\begin{itemize}

\item $[g_0,g_{i,j}]=[g_0,f_i]=[f_s,f_i]=[g_{i,j},g_{l,m}]=e,$ $\forall$ $s,i,l \in \{1,\ldots,k\},\;j,m \in \{1,\ldots,d\}$,

\item $[f_s,g_{i,j}]=g_{i,j-1}^{m_{i,s}}$ $\forall s,i\in \{1,\ldots,k\}$ and $j\in \{2,\ldots,d\}$,

\item $[f_s,g_{i,1}]=g_0^{m_{i,s}}$ $\forall s,i\in \{1,\ldots,k\}.$
\end{itemize}
Note that from the identities $[ab,c]=a[b,c]a^{-1}[a,c]$ and $[a,bc]=[a,b]b[a,c]b^{-1}$, we immediately have the following additional relations
\begin{itemize}
\item $[f_s^{-1},g_{i,j}]\in \langle g_0,g_{i,1},\ldots,g_{i,j-2}\rangle g_{i,j-1}^{-m_{i,s}}\;\forall s,i \in \{1,\ldots,k\}, j\in \{2,\ldots,d\}$,
\item $[f_s^{-1},g_{i,1}]=g_0^{-m_{i,s}}\;\forall s,i \in \{1,\ldots,k\}$.
\end{itemize}

It is easy to see that $G$ is a nilpotent group of degree $d+1$,  and $A=\langle \{g_0\}\cup\{g_{i,j}: (i,j)\in \{1,...,k\}\times\{1,\ldots,d\}\}\rangle$ is a maximal abelian subgroup containing the commutator of $G$ (see Lemma \ref{lem abel max} below). Moreover $k$ is the torsion-free rank of $G/A$, therefore in view of Theorem \ref{teo lower bound} we know that $G$ embeds in $\text{Diff}_+^{1+\alpha}([0,1])$ for $\alpha < 1/k$. To show that $1+1/k$ is actually an upper bound for the regularity we are going to need the following elementary lemma.

\begin{lemma}\label{lem abel max} For all $(i,j)\in \{1,\ldots,k\}\times\{2,\ldots,d\}$ and $n_1,\ldots,n_k\in \mathbb{Z}$ we have 
\begin{enumerate}
\item $[f_1^{n_1}\cdots f_k^{n_k},g_{i,j}]\in \langle g_0,g_{i,1},\ldots,g_{i,j-2}\rangle g_{i,j-1}^{\lambda_i},$
\item $[f_1^{n_1}\cdots f_k^{n_k},g_{i,1}]= g_0^{\lambda_i},$
\end{enumerate}
where $\lambda_i=\sum_{s=1}^k n_s m_{i,s}$. In particular, the subgroup $A$ is a maximal abelian subgroup.

\end{lemma}

\noindent\textbf{Proof.}  To show 1. we do induction on $n=\sum_{s=1}^k \abs{n_s}$. 

Note that when $n=1$ we have the result by the relations of $G$. So, consider an arbitrary natural number $n=\sum_{j=1}^k \abs{n_j}$ and assume that $n_k<0$ (the other case is similar). For all $i \in \{1,\ldots,k\}$ and $j\in\{2,...,d\}$ we have that
$$[f_1^{n_1}\cdots f_k^{n_k},g_{i,j}]=[f_1^{n_1}\cdots f_k^{n_k+1},[f_k^{-1},g_{i,j}]]  \;  [f_k^{-1},g_{i,j}] \; [f_1^{n_1}\cdots f_k^{n_k+1},g_{i,j}],$$
\noindent and since $[f_k^{-1},g_{i,j}]$ belongs to $\langle g_0,g_{i,1},\ldots,g_{i,j-2}\rangle g_{i,j-1}^{-m_{i,k}}$, it follows that $[f_1^{n_1}\cdots f_k^{n_k+1},[f_k^{-1},g_{i,j}]]\in \langle g_0,g_{i,1},\ldots,g_{i,j-2}\rangle$. Also, by induction hypothesis we have $$[f_1^{n_1}\cdots f_k^{n_k+1},g_{i,j}]\in \langle g_0,g_{i,1},\ldots,g_{i,j-2}\rangle g_{i,j-1}^{(\sum_{s=1}^{k-1}n_s m_{i,s}+(n_k+1)m_{i,k})}.$$ Plugging these into the previous equation yields assertion 1. The proof of  assertion 2 is analogous. 

 $\hfill\square$

\begin{remark}  The most useful part of  Lemma \ref{lem abel max} is the explicit expression for the integers $\lambda_i$ appearing. These will be used in the proof of Theorem C.

\end{remark}

\paragraph{Proof of Theorem C.}

\noindent Suppose that $G$ embeds into $\text{Diff}_+^{1+\alpha}([0,1])$ for some $\alpha>1/k$. Let $x_0$ be a point in $(0,1)$ such that $g_{0}(x_0)\neq x_0$ and define the intervals $$I_0:=(\underset{n}{\text{inf}}\; g_0^n(x_0), \underset{n}{\text{sup}}\;g_0^n(x_0))\text{ and }I_{i,j}:=(\underset{n}{\text{inf}}\; g_{i,j}^n(x_0), \underset{n}{\text{sup}}\;g_{i,j}^n(x_0)).$$

\textbf{Case 1:} $f(I_0)\cap I_0=\emptyset$ for all $f\in \langle f_1,\ldots,f_k\rangle\simeq \mathbb{Z}^k$. 

In this case, $I_0$ is a wandering interval for the dynamics of $\langle f_1,\ldots,f_k\rangle$. A contradiction is provided by  Lemma \ref{lem caminatas} followed by Theorem \ref{teo control of distortion} since the central element $g_0$ acts non-trivially on $I_0$.

\vspace{0,2cm}
\textbf{Case 2:} There is  a non-trivial element $f\in \langle f_1,\ldots,f_k\rangle$ such that $f(I_0)= I_0$. 

Let us put $f=f_1^{n_1}\cdots f_k^{n_k}$. Given $i\in\{1,\ldots,k\}$, by Lemma \ref{lem abel max} we have that 
\begin{equation}\label{eq para aplicar conrad}
[f,g_{i,1}]=g_0^{\lambda_i}\text{ and }[f,g_{i,j}]\in \langle g_0,g_{i,1},\ldots,g_{i,j-2}\rangle g_{i,j-1}^{\lambda_i}\hspace{0,3cm}\text{ for all } j\in \{2,\ldots,d\},
\end{equation}
where $\lambda_i=\sum_{j=1}^k n_jm_{i,j}$. Since the vectors $(m_{i,1},\ldots,m_{i,k})$ are linearly independent in $\mathbb{R}^k$, we can choose $i$ to obtain $\lambda_i\neq 0$.
Then, the relations (\ref{eq para aplicar conrad}) and Lemma \ref{lem conrad} imply that $g_{i,1}(I_0)\cap I_0=\emptyset $. Since the action has no crossings,  the element $f$ also fixes the intervals $I_{i,j}$ and hence the same argument also yields that $g_{i,j}(I_{i,j-1})\cap I_{i,j-1}=\emptyset$ for all $j>2$. Therefore, $I_0$ is a wandering interval for the action of $\langle g_{i,1},\ldots, g_{i,k}\rangle\simeq \mathbb{Z}^k$. So, a contradiction is reached using  Lemma \ref{lem caminatas} and Theorem \ref{teo control of distortion}  as before. $\hfill\square$

\subsection{An example with even higher regularity}
\label{sec example with higher regularity}

It is easy to see that in some situations the regularity given by Theorem \ref{teo lower bound} is not critical. In the examples that we know of, this is related to the fact that the group can be split as a direct product of groups each of which allows an embedding with better regularity. Take for example the groups of \cite[\S 4]{int4}. These are given by the presentation
$$G_d:=\langle f,g_1,\ldots ,g_d : [g_i,g_j]=id, [f,g_1]=id, [f,g_{i}]=g_{i-1}\;\forall j\geq 1, i>1 \rangle.$$
Note that $G_d$ is isomorphic to a non-trivial semidirect product of the form $\mathbb{Z}^d\rtimes\mathbb{Z}$.  Now define the group $G:=G_d\times G_d.$  On one hand, it is easy to see that $$G\simeq \mathbb{Z}^{2d}\rtimes \mathbb{Z}^2, $$ and $\mathbb{Z}^{2d}\times \{0\}$ is a maximal abelian subgroup of $G$. Therefore, if we apply Theorem \ref{teo lower bound}, we obtain an embedding of $G$ into $\text{Diff}_+^{1+\alpha}([0,1])$ for all $\alpha< 1/2$. However, on the other hand,  the critical regularity of $G$ is 2. Indeed,  we can apply Theorem \ref{teo lower bound} to each factor of $G$ to obtain an embedding of the factor into  $\text{Diff}_+^{1+\alpha}([0,1])$ for all $\alpha< 1$. If we put these two actions together acting on disjoint intervals (as we did in Section \ref{sec faithful action}), we end up with an embedding of $G$ into $\text{Diff}_+^{1+\alpha}([0,1])$ for all $\alpha< 1$.

\begin{small}

\end{small}

\vspace{0,5cm}

\small 
{Maximiliano Escayola}

IRMAR - UMR CNRS 6625, Université de Rennes, France.

Email: maximiliano.escayola@univ-rennes.fr

\vspace{0,5cm}

Crist\'obal Rivas

Dpto de Matem\'aticas, Universidad de Chile

Las Palmeras 3425, Ñuñoa, Santiago, Chile

Email: cristobalrivas@u.uchile.cl


\begin{thebibliography}{Dillo 83}


\bibitem{bass} {\sc H. Bass.} 
The degree of polynomial growth of finitely generated nilpotent groups.
{\em Proc. Lond. Math. Soc.} {\bf 25} (1972), 603-614.

\bibitem{bergman} {\sc G. Bergman.} 
Right-orderable groups that are not locally indicable. 
{\em Pac. J. Math.} {\bf 147} (1991), 243-248.


\bibitem{BMNR} {\sc C. Bonatti, I. Monteverde, A. Navas \& C. Rivas.} 
Rigidity for $C^1$ actions on the interval arising from hyperbolicity I: solvable groups.
{\em Math. Z.} {\bf 286} (2017), 919-949.


\bibitem{BMRT} {\sc J. Brum, N. Matte Bon, C. Rivas \& M. Triestino.} 
Locally moving groups acting on the line and $\R$-focal actions.
Arxiv:2104.14678 (2021).


\bibitem{calegari} {\sc D. Calegari.} 
Nonsmoothable, locally indicable group actions on the interval. 
{\em Alg. Geom. Top.} {\bf 8} (2008), 609-613.

\bibitem{int4} {\sc G. Castro, E. Jorquera \& A. Navas.} 
Sharp regularity for certain nilpotent group actions on the interval. 
{\em Math. Ann.} {\bf 359} (2014), no. {\bf 1-2}, 101-152.




\bibitem{DKN_acta} {\sc B. Deroin, V. Kleptsyn \& A. Navas.}
Sur la dynamique unidimensionnelle en r\'egularit\'e interm\'ediaire. 
{\em Acta Math.} {\bf 199} (2007), no. {\bf 2}, 199-262.

\bibitem{dknp} {\sc  B. Deroin, V. Klepstyn, A. Navas \& K. Parwani.} Symmetric random walks on $\text{Homeo}_+(\mathbb{R})$.
{\em Ann. Probab.} {\bf 41}, 2066-2089 (May 2013).

\bibitem{GOD} {\sc B. Deroin, A. Navas \& C. Rivas.} Groups, Orders, and Dynamics. ArXiv: 1408.5805.

\bibitem{ff} {\sc B. Farb \& J. Franks.} Groups of homeomorphisms of one-manifolds III: Nilpotent subgroups. {\em Ergod. Th. \& Dynam. Sys,} {\bf 23} (2003), 1467-1484.
  
  
\bibitem{guivarch} {\sc Y. Guivarc'h.} 
Croissance polynomiale et p\'eriodes des fonctions harmoniques. 
{\em Bull. Soc. Math. France} {\bf 101} (1973), 333-379.

\bibitem{jorquera} {\sc E. Jorquera.} A universal nilpotent group of $C^1$ diffeomorphisms of the interval. {\em Topology and its Applications.} {\bf 159} (2012), 2115-2126.
  
\bibitem{JNR} {\sc E. Jorquera, A. Navas \& C. Rivas.} On the sharp regularity for arbitrary actions of nilpotent groups on the interval: the case of $N_4$. {\em Ergod. Th. \& Dynam. Sys.} {\bf 38} (2018), 180-194.




\bibitem{kk} {\sc S. H. Kim \& T. Koberda.} Diffeomorphism groups of critical regularity. {\em Invent. Math.} 221(2) (2020),
421–501.

\bibitem{kk_book} {\sc S. H. Kim \& T. Koberda.} Structure and regularity of group actions on one-manifolds. {\em Springer Monographs in Mathematics} (2021).


\bibitem{kkr} {\sc S. H. Kim ,T. Koberda \&  C. Rivas.} Direct products, overlapping actions and critical regularity. {\em J. Mod. Dyn.} 17 (2021), 285-304.


\bibitem{kopell} {\sc N. Kopell.} Commuting diffeomorphisms. In {\em Global Analisys.} (Berkeley, CA, 1968), Proc. Sympos. Pure Math., Vol XIV, pp 165-184. Amer. Math. Soc. Providence, RI, 1970.


\bibitem{mann-wolff} {\sc K. Mann \& M. Wolff.} 
Reconstructing maps out of groups, \emph{Ann. Sci. \'Ec. Norm. Sup\'er.} {\bf 56}:4 (2023), 1135-1154.


\bibitem{na} {\sc A. Navas.} 
 Groups of Circle Diffeomorphisms. {\em Chicago Lectures in Mathematics.} Univ. of Chicago Press (2011).
 
\bibitem{na1} {\sc A. Navas.} On centralizers of interval diffeomorphisms in critical (intermediate) regularity. {\em J. Anal. Math.} {\bf 121} (2013), 1-30.

\bibitem{navas-loc-ind} {\sc A. Navas.} A finitely generated, locally indicable group with no faithful action by $C^1$ diffeomorphisms of the interval. {\em Geom. \& Topol.} 14, (2010), 573-584.
  

\bibitem{p} {\sc K. Parkhe.} Nilpotent dynamics in dimension one:
Structure and smoothness. {\em  Ergod. Th. \& Dynam. Sys.} {\bf 36} (2016), 2258-2272.

\bibitem{pixton} {\sc D. Pixton.} Nonsmoothable, unstable group actions. {\em Trans. AMS} 229 (1977), 259-268.

\bibitem{pt} {\sc J. Plante \& W. Thurston.} Polynomial growth in holonomy groups of foliations.{\em Comment. Math. Helv.} {\bf 51} (1976),
567-584.


\bibitem{rivas-triestino} {\sc C. Rivas \& M. Triestino.} 
One dimensional actions of Higman's group.
{\em Discrete Analysis} (2019), 15pp.



\bibitem{robinson} {\sc D. J. S. Robinson.} 
A Course in the Theory of Groups. Number 80 in {\em Graduate Texts
in Mathematics.} Springer-Verlag, New York, second edition, 1996.







\bibitem{thurston} {\sc W. Thurston.} 
A generalization of the Reeb stability theorem. 
{\em Topology} {\bf 13} (1974), 347-352.

\bibitem{tsuboi} {\sc T. Tsuboi.} 
Homological and dynamical study on certain groups of Lipschitz homeomorphisms of the circle. 
{\em J. Math. Soc. Japan} {\bf 47} (1995), 1-30.

\end{thebibliography}
\end{document}